\documentclass[a4paper, 12pt, titlepage, fleqn]{article}
\usepackage{times}
\usepackage{amsthm}
\newtheorem{thm}{Theorem}
\newtheorem{cor}{Corollary}

\newtheorem{lemma}[thm]{Lemma}
\newtheorem{prop}{Proposition}

\usepackage{amssymb,amsmath}

\DeclareMathOperator{\F}{\mathbb{F}}

\DeclareMathOperator{\Tr}{Tr}

\begin{document}
		\baselineskip=16.3pt
		\parskip=14pt
		\begin{center}
			\section*{The Number of Zeros of Quadratic Forms with Two Terms}
			{\large 
			Emrah Sercan Y{\i}lmaz \footnote {Research supported by T\"ubitak Grant 117F274} 
			\\ 
			Bo\u gazi\c ci University,
			Turkey}

		\end{center}

	\subsection*{Abstract}
	In this paper, we will interest in finding the number of zeros of the quadratic forms over finite fields. We will apply the tool for finding the number of rational points of supersingular curves in \cite{MY-reduction}. We will give some more tools for Artin-Schreier curves related with quadratic form.  We will especially interest in the (supersingular) Artin-Schreier curve $y^q-y=x^{q^b+1}-x^{q^a+1}$.

\textbf{Keywords:} Number of Rational Points, Supersingular Curves, Quadratic Forms.

\section{Introduction}
Quadratic forms over finite fields are using in coding theory and also appear in proofs of many interesting problems.  The number of rational points of curves  $y^q-y=x^{q^b+1}-x^{q^a+1}$ over $\F_{q}$ are studied in a paper by Lahtonen et al. for $q=2$ with some conditions on $a$ and $b$  and are presented in a paper by Fitzgerald without any condition. They are presented in a paper by Roy for $q=2^r$ where $r$ is a positive integer. The curves $y^q-y=x^{q+1}-x^2$ are related to the number of irreducible polynomials with prescribed coefficients over finite fields by Ahmadi et al for even $q$ with some other quadratic forms and  by McGuire-Y\i lmaz for odd $q$.  In this paper, we will give some general tools for finding the zeros of quadratic forms over finite fields and we will especially interest in the curves $y^q-y=x^{q^b+1}-x^{q^a+1}$ over $\F_q$ where $q$ is an odd prime power and $a$, $b$ are positive integer having different $2$-adic and $p$-adic valuation. Note that the methods in this paper will also allow to find the number of zeros of more complicated quadratic forms over finite fields.

Let $q$ be a prime $p$ power and $n$, $a$ and $b$ be positive integers in this paper. We denote the greatest common divisor of $a$ and $b$ with $(a,b)$ and the least common multiple of $a$ and $b$ with $[a,b]$. We denote prime $\ell$-adic valuation of $a$ with $\nu_\ell(a)$.
\section{Preliminaries}

\subsection{Supersingular Curves}
A protective smooth absolutely irreducible curve $X$ of genus $g$ defined over $\mathbb F_q$ is \emph{supersingular} if all Weil numbers of $X$ have the form $\eta_i = \sqrt{q}\cdot \zeta_i$ where $\zeta_i$ is a root of unity. For any $n\ge 1$ we have \begin{equation}\label{eqn-root-sum-2}
-q^{-n/2}[\#X(\F_{q^n})-(q^n+1)]=\sum_{i=1}^{2g} \zeta_i^{n}.
\end{equation}

The smallest positive integer $s=s_X$ such that $\zeta_i^{s}=1$ for all 
$i=1,\ldots	,2g$ will be called the \emph{period} of $X$. 
The period depends on $q$, in the sense that $X(\F_{q^n})$
may have a different period to $X(\F_{q})$.

If all $\zeta_i^{n}$ are $1$, we call the curve $X$ minimal on $\F_{q^n}$. 
If all $\zeta_i^{n}$ are $-1$, we call the curve $X$ maximal on $\F_{q^n}$. Every supersingular curve are minimal on some extensions of $\F_q$. On the other hand, every curve which attains maximal or minimal on an extension of $\F_q$ is supersingular.
\subsection{Basics on Quadratic Forms}

Let $L(x)$ be a $\F_q$-linearized polynomial. It is known that 
the functions 
from $\F_{q^n}$ to $\F_q$ 
such that
$
x\mapsto \Tr_{\F_{q^n}/\F_q}\left(xL(x)\right)
$
is a quadratic form over $\F_q$.
The number of zeros of such a quadratic form $Q$
can be written as 
$
q^{n-1}+\lambda(q-1)q^{\frac{n+w}{2}-1}
$
where $\lambda\in \{-1,0,1\}$
and $w$ is the dimension of the radical 
$\{x\in \F_{q^n} \ : \ Q(x+y)-Q(x)-Q(y)=0 \text{ for all } y \in \F_{q^n} \}$
of $Q$ when $q$ is odd.
The dimension $w$ has a slightly different definition for even characteristic. We assume that the characteristic is odd.  
Note that if $n$ and $w$ have different parity (odd/even), then $\lambda$ has to be $0$. Otherwise, it is well-known that $\lambda$ is nonzero.

The quadratic form 
$\Tr_{\F_{q^n}/\F_q}\left(xL(x)\right)$
is related with the curve 
$
C_L:y^q-y=xL(x).
$
The number of $\F_{q^n}$-rational points of this curve is $qN+1$ 
where $N=\{x\in \F_{q^n} \ : \ \Tr_{\F_{q^n}/\F_q}\left(xL(x)\right)=0\}$.
We will write the number of $\F_{q^n}$-rational points of $C_L$ as
$$\#C_L(\F_{q^n})=q^n+1+\lambda_n(C_L)(q-1)q^{\frac{n+w_n(C_L)}{2}}$$
where $\lambda_n(C_L)\in \{-1,0,1\}$
and $w_n(C_L)$ is the dimension of the radical
of $\Tr_{\F_{q^n}/\F_q}\left(xL(x)\right)$. We note that the curves related with quadratic forms are supersingular curves. Fore more details, see \cite{Van} or  \cite{MY-slope}.

\section{Some Reduction Theorems for Supersingular Curves}
The following theorem was proven in \cite{MY-reduction} and it is an effective tool to find the number of rational points of supersingular curve over finite fields. 
\begin{thm}\label{reduction-thm}
	Let $X$ be a supersingular curve of genus $g$ defined over $\mathbb F_q$ with
	period $s$.
	Let $n$ be a positive integer, let $\gcd (n,s)=m$ and write $n=m\cdot t$. If $q$ is odd, then we have $$
	\#X(\F_{q^n})-(q^n+1)=\begin{cases}
	q^{(n-m)/2}[\#X(\F_{q^m})-(q^m+1)] &\text{if } m\cdot r \text{ is even},	\\
	q^{(n-m)/2}[\#X(\F_{q^m})-(q^m+1)]\left(\frac{(-1)^{(t-1)/2}t}{p}\right)&\text{if } m \cdot r \text{ is odd and } p\nmid t,
	\\
	q^{(n-m)/2}[\#X(\F_{q^m})-(q^m+1)]&\text{if } m \cdot r \text{ is odd and } p\mid t.
	\end{cases}$$ If $q$ is even, then we have  $$
	\#X(\F_{q^n})-(q^n+1)=\begin{cases}
	q^{(n-m)/2}[\#X(\F_{q^m})-(q^m+1)] &\text{if } m\cdot r \text{ is even},\\
	q^{(n-m)/2}[\#X(\F_{q^m})-(q^m+1)](-1)^{(t^2-1)/8}&\text{if } m \cdot r \text{ is odd}.\\
	\end{cases}$$
\end{thm}
This theorem can be easily generalised for maximal curves. Since manipulating the related proof in \cite{MY-reduction} is straightforward, we will skip the proof. The proof of Theorem \ref{reduction-thm-m} depends of the fact that $r$-th and $2r$-th cyclotomic polynomials $\Phi_r$ and $\Phi_{2r}$ have the relation $\Phi_{2r}(x)=\Phi_r(-x)$ where $r$ is odd. 
\begin{thm}\label{reduction-thm-m}
	Let $X$ be a  maximal curve of genus $g$ defined over $\mathbb F_q$ with period $s=2s^\prime$.
	Let $n$ be a positive integer with $\gcd (n,s^\prime) \ne \gcd (n,s)$. Then we have $$
	\#X(\F_{q^n})-(q^n+1)=	-q^{n/2}[\#X(\F_{q^{n/2}})-(q^{n/2}+1)].$$
\end{thm}

There is also a straightforward result for maximal curves. See \cite{MY-L divisibility}. We will use contrapositive of this theorem to show that a curve is not maximal on a given extension when we are try to distinguish that it is maximal or minimal.

\begin{prop}
	Let $X$ be a curve defined over $\F_q$. If $X(\F_{q^{2n}})$ is maximal, then $\#X(\F_{q^n})=q^n+1$. 
\end{prop}

\begin{proof}
	If $X(\F_{q^{2n}})$ is maximal then each Weil number of $X$ over $\mathbb F_{q^{2n}}$ is $-1$. Therefore,  Weil number of $X$ over $\mathbb F_{q^{n}}$ is either $i$ or $-i$. A $\mathbb Z$ linear combination of the set $\{\pm i\}$ is a real number only if it is $0$.    
\end{proof}

\section{Some Reduction Theorems for Artin-Schreier Curves related with Quadratic Forms}
In this section we will give some additional tools to find zeros of quadratic forms over finite fields in addition to Theorem \ref{reduction-thm} and \ref{reduction-thm-m} for odd characteristics.

The following proposition is well-known. See \cite{Lidl} for example.
\begin{prop}\label{Q=c}
	Let $L(x)$ be a $\F_q$-linearized polynomial and define $Q(x)=\Tr_{\F_{q^n}/\F_q}(xL(x))$. Assume that the set $\{x\in \F_{q^n} \ : \ Q(x)=0 \}$ has $q^{n-1}+\lambda(q-1)q^{\frac{n+w}{2}-1}$ elements where $w$ is the dimension of radical of $Q(x)$ with $\lambda\in\{-1,+1\}$. Then the set  $\{x\in \F_{q^n} \ : \ Q(x)=c \}$ has $$q^{n-1}-\lambda q^{\frac{n+w}{2}-1}$$ elements  for any nonzero $c\in \F_q$.
\end{prop}

\begin{thm}
	Let $Q=q^d$ with some positive integer $d$ and $L(x)$ be a $\F_Q$-linearized polynomial. If the curve $y^Q-y=xL(x)$ has a nonzero sign $\lambda$ over $\F_{Q^n}$, then  $y^q-y=xL(x)$ has  the same sign $\lambda$ over $\F_{q^{dn}}$.
\end{thm}

\begin{proof}
	This result is a consequence of the set equality $$\{ x\in \F_{q^{dn}} \ : \ \Tr_{\F_{q^n}/\F_q}(xL(x))=0 \}=\bigcup\limits_{\begin{matrix}
		c\in \F_Q \text{ with} \\ 	\Tr_{\F_{Q}/\F_q}(c)=0
		\end{matrix}} \{ x\in \F_{q^{dn}} \ : \ \Tr_{\F_{q^n}/\F_Q}(xL(x))=c \}.$$	Let the dimension of radical of $L(x)$ over $\F_Q$ be $n$ and let the sign of  $y^q-y=xL(x)$ be $\lambda^\prime$. By Proposition \ref{Q=c} and the set equality above we have 
	$$
	q^{nd-1}+\lambda^\prime(q-1)q^{\frac{nd+wd}{2}-1} = (Q^{n-1}+\lambda(Q-1)Q^{\frac{n+w}{2}-1})+\left(\frac Qq-1\right)(Q^{n-1}-\lambda Q^{\frac{n+w}{2}-1}).$$ Hence, $\lambda^\prime=\lambda$.
\end{proof}

\begin{thm}\label{Qq}
	Let $\ell$ be a prime with $(\ell,2p)=1$. Let $C:y^q-y=xL(x)$ with $\F_q$-linearized polynomial $L(x)$.  If $\lambda_{n}(C),\lambda_{n\ell}(C)\ne 0$, then we have the equality $$\lambda_{n}(C)\equiv \lambda_{n\ell}(C) q^{\frac{n(l-1)+w_{nl}(C)-w_n(C)}{2} }\mod \ell.$$
\end{thm}

\begin{proof}
	Let $\alpha\in \mathbb F_{q^n}$ with $\Tr_{\F_{q^n}/\F_q}(\alpha)=c\ne 0$. Define $C_\alpha:y^q-y=xL(x)-\alpha$ over $\F_{q^{n}}$. Let $\beta \in C_\alpha(\F_{q^{n\ell}})-C_\alpha(\F_{q^n})$. Then $\beta,\beta^q,\ldots,\beta^{q^{n(\ell-1)}}$ are also in $ C_\alpha(\F_{q^{n\ell}})-C_\alpha(\F_{q^n})$. Therefore, $\ell$ divides $\#C_\alpha(\F_{q^{n\ell}})-\#C_\alpha(\F_{q^n})$. By Proposition \ref{Q=c} we have $$\#C(\mathbb F_{q^n})=q^n+1-\lambda_{n}(C)q^{\frac{n+w_n(C)}{2}} \ \ \ \text{ and } \ \ \ \#C(\mathbb F_{q^{n\ell}})=q^n+1-\lambda_{n\ell}(C)q^{\frac{n\ell+w_{n\ell}(C)}{2}}.$$ This yields the result.
\end{proof}
Note:  Our reason to interest in the curve $C_\alpha$ is to get rid of the multiplicity $q-1$. The same proof works for the curve $C$ and the prime $\ell$ with an additional condition that $\ell$ does not divide $q-1$.

\begin{cor}\label{cor}
	Let $n$ be a positive even integer and let $t$ be a positive integer with $(t,2p)=1$. Let $C:y^q-y=xL(x)$ with $\F_q$-linearized polynomial $L(x)$. Assume that $w_n(C)$ is even and $w_{nd}(C)=d w_{n}(C)$ for all $d\mid t$. Then $\lambda_{nd}(C)=\lambda_{n}(C)$ for all $d\mid t$. 
\end{cor}
\begin{proof}
	Let $\ell$ be a prime with $d\mid t$ and  $d\ell\mid t$. By Theorem \ref{Qq} we have that $$\lambda_{nd}(C) \equiv q^{\frac{nd+w_{nd}(C)}{2}(l-1)}\lambda_{nd\ell}\equiv \lambda_{nd\ell}\mod \ell .$$ Therefore, we have  $\lambda_{nd}(C)=\lambda_{nd\ell}(C).$ This finishes the proof.
\end{proof}
\section{Main Result}
%Since $(x^{b+a}-1)(x^{b-a}-1)$ divides $x^{pl}-1$ where $l=[b+a,b-a]$, the ... is maximal.

In this section, we mainly prove the theorem \ref{thm-ba}. This theorem is depending on Theorems \ref{thm-b0-0}, \ref{thm-b0-1} and \ref{thm-b0-2} which are related with the number of rational points of the curve $y^q-y=x^{q^b+1}-x^2$ over $\F_q$ for any positive integer $b$.

In order to find these numbers we need to find dimension of radicals and related signs. Finding the dimension of radicals are straightforward for the class we are interesting. We will find them in the following lemma. Therefore, the problem will be reduced to find the related signs. 

\begin{lemma}\label{lem1}
	The dimension of the radical of $\Tr_{\F_{q^n}/\F_q}(x^{q^b+1}-x^{q^a+1})$ is $$\begin{cases}
		(b+a,n)+(|b-a|,n)-(b+a,|b-a|,n) &\text{if } \nu_p(n) \le  \max\{\nu_p(b+a),\nu_p(|b-a|)\},\\
		(b+a,n)+(|b-a|,n) &\text{if } \nu_p(n)  > \max\{\nu_p(b+a),\nu_p(|b-a|)\}.
	\end{cases}$$
\end{lemma}

\begin{proof}
	Without loss of generality we will assume that $b>a$. Then the radical is the set \begin{align*}
	&\{x\in \F_{q^n} \ : \ \Tr_{\F_{q^n}/\F_q}(x^{q^b}y-x^{q^a}y+y^{q^b}x-y^{q^a}x)=0 \text{ for all } y \in \F_{q^n} \}\\
	=&\{x\in \F_{q^n} \ : \ \Tr_{\F_{q^n}/\F_q}(y^{q^b}(x^{q^{2b}}-x^{q^{a+b}}+x-x^{q^{b-a}}))=0 \text{ for all } y \in \F_{q^n} \}\\
	=&\{x\in \F_{q^n} \ : \ x^{q^{2b}}-x^{q^{a+b}}-x^{q^{b-a}}+x=0 \}.
	\end{align*}
	Therefore the dimension of the radical is $$(x^{2b}-x^{b+a}-x^{b-a}+1,x^n-1)=((x^{b+a-1})(x^{b-a}-1),x^n-1).$$The result follows by this equality.
\end{proof}

The following lemma was proved in \cite{MY-irred-p-2term}. 
\begin{lemma}\label{lem2}
	The curve $y^q-y=x^{q+1}-x^2$ over $\F_{q^{2pt}}$ is minimal if $q\equiv 1 \mod 4$ and maximal if $q\equiv 3 \mod 4$.
\end{lemma}
The idea  of  proof was depending on quadratic Gauss sum. 
Same reasoning also gives the following result. 
\begin{lemma}
	The curve $y^q-y=x^{q^2+1}-x^2$ over $\F_{q^{2pt}}$ is minimal if $q\equiv 1 \mod 4$ and maximal if $q\equiv 3 \mod 4$.
\end{lemma}

\begin{thm}\label{thm-b0-0}
	Let $b$ be a positive odd integer with $\nu_p(b)=\ell$.  Let $n$ be positive integer with $(n,4pb)=m$ and write $n=m\cdot t$. Let $C_{q,b,0}:y^q-y=x^{q^b+1}-x^2$ be a curve defined over $\F_{q}$. Then the following equalities satisfy.
	$$\lambda_{n}(C_{q,b,0})=\begin{cases}
	\left(\frac{(-1)^{(t-1)/2}t}{p}\right)& \text{if }\nu_2(m)=0 \text{ and } \nu_p(m)\le \ell,\\
	0& \text{if } \nu_2(m)=0 \text{ and } \nu_p(m)=\ell+1,\\
	0& \text{if } \nu_2(m)=1 \text{ and } \nu_p(m)\le \ell,\\
	(-1)^{\frac{q+1}2}&\text{if } \nu_2(m)= 1 \text{ and } \nu_p(m)=\ell+1.\\
	0&\text{if } \nu_2(m)= 2 \text{ and } \nu_p(m)\le \ell,\\
	-1&\text{if } \nu_2(m)= 2 \text{ and } \nu_p(m)=\ell+1.
	\end{cases}$$	
\end{thm}
%\left(\frac{(-1)^{(c-1)/2}c}{p}\right)
\begin{proof} The degree of $((x^b-1)^2,x^{2pb}-1)$ is $2b$. Since $\F_q$-dimension of the radical of  $C_{q,b,0}$ over $\F_{q^{2pb}}$ is maximal and  $\F_q$-codimension of the radical  of  $C_{q,b,0}$ over $\F_{q^{2pb}}$ is even, the curve $C_{q,b,0}$ is maximal or minimal over $\F_{q^{2pb}}$. We will prove the theorem for $m$ dividing $2pd$ and then use Theorem \ref{reduction-thm} and \ref{reduction-thm-m} to complete for all positive integer $n$.
	 
	\underline{Case $\nu_2(m)=0$ and  $\nu_p(m)\le \ell$:} \\
	In this case, $m$ divides $b$ and hence $x^{q^b+1}=x^2$ for all $x\in \F_{q^m}$. Therefore, $$\#X(\F_{q^n})=q\cdot q^n+1=q^d+1+(q-1)q^{\frac n2}q^{\frac n2}.$$
	
	\underline{Case $\nu_2(m)=0$ and  $\nu_p(m)=\ell+1$:} \\
	We have $m-\deg((x^b-1)^2,x^{m}-1)=m-2(m/p)=(p-2)\cdot m/p$ is odd. 
	
	\underline{Case $\nu_2(m)=1$ and  $\nu_p(m)\le \ell $:} \\
	We have $m-\deg((x^b-1)^2,x^{m}-1)=m-(m/2)=m/2$ is odd.
	
	\underline{Case $\nu_2(m)=1$ and  $\nu_p(m)=\ell+1$:} \\
	Since $\deg((x^b-1)^2,x^{md}-1)=d\deg((x^b-1)^2,x^{m}-1)$ for all $d\mid (2pb/m)$, $y^{q}-y=x^{ q^b+1}-x^2$ over $\F_{q^{m}}$ and $\F_{q^{2pb}}$ has same sign by Corollary \ref{cor}.
	
	Moreover, $y^{q}-y=x^{ q^b+1}-x^2$ over $\F_{q^{2pb}}$ and  $y^{q^b}-y=x^{ q^b+1}-x^2$ over $\F_{(q^{b})^{2p}}$ has same sign by Theorem \ref{Qq} since both has nonzero sign. 
	
	Indeed, $y^{q^b}-y=x^{ q^b+1}-x^2$ over $\F_{(q^{b})^{2p}}$ has sign $	(-1)^{\frac{q^b+1}2}=(-1)^{\frac{q+1}2}$ by Lemma \ref{lem1}.
\end{proof}

\begin{thm}\label{thm-b0-1}
		Let $b$ be a positive odd integer with $\nu_p(b)=\ell$.  Let $n$ be positive integer  with $(n,4pb)=m$ and write $n=m\cdot t$. Let $C_{q,2b,0}:y^q-y=x^{q^{2b}+1}-x^2$ be a curve defined over $\F_{q}$. Then the following equalities satisfy.
	$$\lambda_{n}(C_{q,2b,0})=\begin{cases}
	\left(\frac{(-1)^{(t-1)/2}t}{p}\right)& \text{if }\nu_2(m)=0 \text{ and } \nu_p(m)\le \ell,\\
	0& \text{if } \nu_2(m)=0 \text{ and } \nu_p(m)= \ell+1,\\
	1& \text{if } \nu_2(m)=1 \text{ and } \nu_p(m)le  \ell,\\
	(-1)^{\frac{q+1}2}&\text{if } \nu_2(m)= 1 \text{ and } \nu_p(m)= \ell+1.\\
	(-1)^{\frac{q-1}{2}}& \text{if } \nu_2(m)=2 \text{ and } \nu_p(m)\le \ell,\\
	-1&\text{if } \nu_2(m)= 2 \text{ and } \nu_p(m)= \ell+1.
	\end{cases}$$	
\end{thm}

\begin{proof} The degree of $((x^{2b}-1)^2,x^{2pb}-1)$ is $4b$. Since $\F_q$-dimension of the radical of  $C_{q,2b,0}$ over $\F_{q^{2pb}}$ is maximal and  $\F_q$-codimension of the radical  of  $C_{q,2b,0}$ over $\F_{q^{2pb}}$ is even, the curve $C_{q,2b,0}$ is maximal or minimal over $\F_{q^{2pb}}$. We will prove the theorem for $m$ dividing $2pd$ and then use Theorem \ref{reduction-thm} and \ref{reduction-thm-m} to complete for all positive integer $n$.
	
	\underline{Case $\nu_2(m)\in\{0,1\}$ and  $\nu_p(m)\le \ell$:} \\
	In this case, $m$ divides $b$ and hence $x^{q^{2b}+1}=x^2$ for all $x\in \F_{q^m}$. Therefore, $$\#X(\F_{q^n})=q\cdot q^n+1=q^d+1+(q-1)q^{\frac n2}q^{\frac n2}.$$
	
	\underline{Case $\nu_2(m)=0$ and  $\nu_p(m)= \ell+1$:} \\
	We have $m-\deg((x^{2b}-1)^2,x^{m}-1)=m-2(m/p)=(p-2)\cdot m/p$ is odd.
	
	\underline{Case $\nu_2(m)=1$ and  $\nu_p(m)= \ell+1$:} \\
	Since $\deg((x^{2b}-1)^2,x^{md}-1)=d\deg((x^{2b}-1)^2,x^{m}-1)$ for all $d\mid (2pb/m)$, $y^{q}-y=x^{ q^{2b}+1}-x^2$ over $\F_{q^{m}}$ and $\F_{q^{2pb}}$ has same sign by Corollary \ref{cor}.
	
	Moreover, $y^{q}-y=x^{ q^{2b}+1}-x^2$ over $\F_{q^{2pb}}$ and  $y^{q^{b}}-y=x^{ q^{2b}+1}-x^2$ over $\F_{(q^{b})^{2p}}$ has same sign by Theorem \ref{Qq} since both has nonzero sign. 
	
	Indeed, $y^{q^b}-y=x^{ q^{2b}+1}-x^2$ over $\F_{(q^{b})^{2p}}$ has sign $	(-1)^{\frac{q^b+1}2}=(-1)^{\frac{q+1}2}$ by Lemma \ref{lem2}.
\end{proof}

\begin{thm}\label{thm-b0-2}
	Let $b$ be a positive odd integer with $\nu_p(b)=\ell$.  Let $n$ and $k\ge 2$ be positive integers with $(n,2^kpb)=m$ and write $n=m\cdot t$. Let $C_{q,2^kb,0}:y^q-y=x^{q^{2^kb}+1}-x^2$ be a curve defined over $\F_{q}$. Then the following equalities satisfy.
	$$\lambda_{n}(C_{q,2^kb,0})=\begin{cases}
	\left(\frac{(-1)^{(t-1)/2}t}{p}\right)& \text{if }\nu_2(m)=0 \text{ and } \nu_p(m)\le \ell,\\
	0& \text{if } \nu_2(m)=0 \text{ and } \nu_p(m)= \ell+1,\\
	1& \text{if } \nu_2(m)=1 \text{ and } \nu_p(m)\le \ell,\\
	(-1)^{\frac{q+1}2}&\text{if } \nu_2(m)= 1 \text{ and } \nu_p(m)=\ell+1.\\
	(-1)^{\frac{q-1}{2}}& \text{if } \nu_2(m)=2 \text{ and } \nu_p(m)\le \ell,\\
	-1&\text{if } \nu_2(m)= 2 \text{ and } \nu_p(m)=\ell+1.
	\end{cases}$$	
\end{thm}

\begin{proof} The degree of $((x^{2^kb}-1)^2,x^{2^kpb}-1)$ is $2^{k+1}b$. Since $\F_q$-dimension of the radical of  $C_{q,2^kb,0}$ over $\F_{q^{2pb}}$ is maximal and  $\F_q$-codimension of the radical  of  $C_{q,2b,0}$ over $\F_{q^{2^kpb}}$ and $\F_{q^{2^{k-1}pb}}$   are even, the curve $C_{q,2^kb,0}$ is minimal over $\F_{q^{2pb}}$. We will prove the theorem for $m$ dividing $2pd$ and then use Theorem \ref{reduction-thm} to complete for all positive integer $n$.
	
	\underline{Case $\nu_2(m) \le k$ and  $\nu_p(m)\le \ell$:} \\
	In this case, $m$ divides $b$ and hence $x^{q^{2^kb}+1}=x^2$ for all $x\in \F_{q^m}$. Therefore, $$\#X(\F_{q^m})=q\cdot q^m+1=q^d+1+(q-1)q^{\frac m2}q^{\frac m2}.$$
	
	\underline{Case $\nu_2(m)=0$ and  $\nu_p(m)=\ell+1$:} \\
	We have $m-\deg((x^{2b}-1)^2,x^{m}-1)=m-2(m/p)=(p-2)\cdot m/p$ is odd.
	
	\underline{Case $\nu_2(m)\ge 2$ and  $\nu_p(m)=\ell+1$:} \\
	Since $\deg((x^{2^kb}-1)^2,x^{md}-1)=d\deg((x^{2^kb}-1)^2,x^{m}-1)$ for all $d\mid (2^{\nu_2(m)}pb/m)$, $y^{q}-y=x^{ q^{2^kb}+1}-x^2$ over $\F_{q^{m}}$ and $\F_{q^{2^{\nu_2(m)}pb}}$ has same sign by Corollary \ref{cor} and has same sign with  $\F_{q^{2^{k}pb}}$ by Theorem \ref{Qq} which is $-1$.  
	
	\underline{Case $\nu_2(m)=1$ and  $\nu_p(m)=\ell+1$:} \\
  $y^{q}-y=x^{ q^{2^kb}+1}-x^2=x^{ q^{2p(2^{k-1}b+1)}+1}-x^2$ over $\F_{q^{2p}}$ has sign $(-1)^{p+1}$   by Theorem \ref{thm-b0-1}.
	
\end{proof}

\begin{thm}\label{thm-ba}
	Let $a$ and $b$ be a positive integers with $\nu_2(b)<\nu_2(a)=\ell$ and $\nu_2(a)\ne \nu_2(b)=k$ and  let $n$ be positive integers with $(n,s)=m$ where $s$ is $2p[b+a,|b-a|]$ if $\ell=1$ and $p\equiv 1\mod 4$ and $p[b+a,|b-a|]$ otherwise.  Write $n=m\cdot t$. Let $C_{q,b,a}:y^q-y=x^{q^{b}+1}-x^{q^a+1}$ be a curve defined over $\F_{q}$. Then the following equalities satisfy.
	$$\lambda_{n}(C_{q,b,a})=\begin{cases}
	\left(\frac{(-1)^{(t-1)/2}t}{p}\right)& \text{if }\nu_2(m)=0 \text{ and } \nu_p(m)\le \ell,\\
	1& \text{if }\nu_2(m)\ge 1 \text{ and } \nu_p(m)\le \ell,\\
	0& \text{if } \nu_2(m)=0 \text{ and } \nu_p(m)=\ell+1,\\
	\lambda_{n}(C_{q,b,0})& \text{if } \nu_2(m)\ge 1 \text{ and } \nu_p(m)=\ell+1.\\
	\end{cases}$$	
\end{thm}

\begin{proof}
	\underline{Case $\nu_2(m) \le k$ and  $\nu_p(m)\le \ell$:} \\
In this case, $m$ divides $b$ and hence $x^{q^b+1}-x^{p^a+1}=0$ for all $x\in \F_{q^m}$. Therefore, $$\#C_{q,b,a}(\F_{q^n})=q\cdot q^n+1=q^d+1+(q-1)q^{\frac n2}q^{\frac n2}.$$

\underline{Case $\nu_2(m)=0$ and  $\nu_p(m)=\ell+1$:} \\
We have $m-\deg((x^{b+a}-1)(x^{|b-a|}-1),x^{m}-1)=m-[(m,b+a)+(m,|b-a|)$ is odd.

\underline{Case $\nu_2(m)\ge 1$ and  $\nu_p(m)=\ell+1$:} \\
In this case, $m$ divides $a$ and hence $x^{q^a+1}=x^2$ for all $x\in \F_{q^m}$. Therefore, $$\#C_{q,b,a}(\F_{q^m})=\#C_{q,b,0}(\F_{q^m}).$$
	
\end{proof}

\end{document}